\documentclass{article}

\usepackage{a4wide}
\usepackage{xypic}
\usepackage{amsmath}
\usepackage{amsfonts}
\usepackage{amssymb}
\usepackage{amsthm}
\usepackage{bbm}

\newtheorem{theorem}{Theorem}[section]

\newtheorem{lemma}[theorem]{Lemma}
\newtheorem{proposition}[theorem]{Proposition}
\newtheorem{corollary}[theorem]{Corollary}

\newtheorem{example}[theorem]{Example}
\newtheorem{definition}[theorem]{Definition}

\newcommand{\xf}{\mathbb{F}}
\newcommand{\xg}{\mathbb{G}}
\newcommand{\xq}{\mathbb{Q}}
\newcommand{\xp}{\mathbb{P}}
\newcommand{\xa}{\mathbb{A}}
\newcommand{\xh}{\mathbb{H}}
\newcommand{\xz}{\mathbb{Z}}
\newcommand{\pol}{\mathcal{L}}
\newcommand{\emm}{\mathcal{M}}
\newcommand{\inzweizweimat}[4]{\begin{array}{c@{\quad}c}
#1 & #2 \\
#3 & #4
\end{array}}
\newcommand{\squaremat}[4]{\left( \inzweizweimat{#1}{#2}{#3}{#4} \right)}

\title{Theta relations with real multiplication \\
by square root three}

\author{Robert Carls}

\begin{document}

\maketitle

{\small 
\begin{center}
\begin{tabular}{l@{\hspace{0.5cm}}l}
Contact details: \\
Universit\"{a}t Ulm,
Institut f\"{u}r Reine Mathematik \\
D-89069 Ulm, Germany, \texttt{robert.carls@uni-ulm.de} \\
Ph.: +49(0)7315023564, Fax: +49(0)7315023618
\end{tabular}
\end{center}
}

\begin{abstract}
\noindent
In this article we compute new theta relations which define
the moduli space of abelian surfaces with real multiplication
by square root three.
We give the locus of square root three abelian surfaces
in terms of the canonical coordinates on
Mumford's moduli space of abelian surfaces with theta
structure.
The method that we use to compute the real
multiplication theta relations is of purely
algebraic nature. It differs from Runge's method and
recent work by Gruenewald. We expect that our method
generalizes to higher discriminants.
\end{abstract}

\noindent
{\small
\begin{center}
\begin{tabular}{ll}
KEYWORDS: & Real Multiplication, Theta Functions, Hilbert-Blumenthal-Varieties \\
AMS MSC: & 11F41, 11G10
\end{tabular}
\end{center}
}

\section{Introduction}
\label{intro}

The subject of real multiplication is a classical one.
Many articles and books (see for example \cite{vg88} \cite{fr90})
have been written on this subject.
In the center of interest are the so called
Hilbert-Blumenthal varieties which were introduced by Hilbert
and subsequently studied by Hilbert's
student Blumenthal in his Habilitationsschrift.
The Hilbert-Blumenthal varieties form classifying spaces for abelian surfaces
whose ring of endomorphisms contains a specific order of a given real quadratic field. Their images in the moduli space of
abelian surfaces are sometimes called Humbert surfaces.
The classical theory uses the complex analytic
description of the latter moduli spaces as quotients of
a generalized upper half plane by the action of the Hilbert
modular group.
\newline\indent
The arithmetic of the `special points' on Hilbert-Blumenthal varieties
is completely determined by the corresponding Shimura data.
In the case of Hilbert-Blumenthal varieties a suitable Shimura data $(G,X)$ can
be made up as follows .
Let $K$ be a real quadratic number field and denote by $O_K$
its ring of integers. We may assume that $K=\xq(\sqrt{d})$, where
$d \geq 3$ is an odd square-free integer.
The reductive group $G$ is given by the pull back
along the morphisms
$\mathrm{Res}_{O_K/\xz} \mathrm{GL}_{2,O_K}
\stackrel{\mathrm{det}}{\rightarrow} \mathrm{Res}_{O_K/\xz} \xg_{m,O_K}$
and $\xg_{m,\xz} \hookrightarrow \mathrm{Res}_{O_K/\xz} \xg_{m,O_K}$.
The space $X$ is given by the disjoint union of the
double upper and lower half plane $(\xh^{+})^2 \cup (\xh^{-})^2$.
We define a compact open subgroup $H$ of
$G(\hat{\xz})$ by setting $H
= \mathrm{Ker} \big[ G(\hat{\xz}) \rightarrow G(\xz / 4 \xz) \big]$.
It follows from Deligne's theory that
there exists an algebraic variety $\mathrm{Sh}(G,X)_{H}$ such that
its complex points are parametrized by the set
$G(\xq) \backslash X \times G \big( \xa_{\xq}^f \big) / H$,
where $\xa_{\xq}^f$ denotes the finite adeles of $\xq$,
the group $G(\xq)$ acts on $X$ by fractional linear transformations and
by translation on the second factor.
The points of $\mathrm{Sh}(G,X)_{H}$
correspond to quadruples $(A,\alpha,\lambda,\eta)$, where
\begin{enumerate}
\item
$A$ is an abelian variety,
\item
$\alpha:O_K \rightarrow \mathrm{End}(A)$ is a morphism
of rings,
\item
$\lambda:A \rightarrow \mathrm{Pic}^0_A$ is a principal
$O_K$-linear polarization,
\item
$\eta: (O_K/4 O_K)^2 \stackrel{\sim}{\rightarrow}
A[4]$ is an isomorphism of $O_K$-modules that
respects natural pairings.
\end{enumerate}
Our aim is to give algebraic equations over $\xz$
which define a somewhat modified functorial version of
the Shimura variety $\mathrm{Sh}(G,X)_{H}$.
The special case of $K=\xq(\sqrt{3})$ is discussed
below.
We note that the $O_K$-linearity of the polarization
$\alpha$ will reappear in the algebraic context as
$\sqrt{d}$-admissibility of certain line bundles.
The compatibility of the level structure
$\eta$ with the $O_K$-action, which implies that the action
is diagonalized, is built in our
proposed notion of $\sqrt{d}$-admissible theta structure.
\newline\indent
We study the subject of real multiplication from
an algebraic point of view.
Our method is of purely algebraic
nature, based on the combinatorics of algebraic
theta functions. The concept of
algebraic theta functions was introduced by Mumford
in \cite{mu66} \cite{mu67a} \cite{mu67b}.
Mumford writes in the
introduction of his article \cite{mu66} the following
`There are several interesting topics which I have not gone into in this paper, but which can be investigated in the same spirit: for example ...
an analysis of special theta-functions and special abelian varieties; ...'.
Among others, one would like to
have an explicit description of the moduli space
of abelian varieties with a given multiplication
type in terms of Mumford's canonical coordinate system
on the moduli space of abelian varieties with theta structure.
In this article we solve this problem in the special
case of abelian surfaces with $\sqrt{3}$-multiplication
and level-$4$ theta structure.
We also treat the `degenerate' case, namely the case
of split abelian surfaces.
\newline\indent
Regarding the problem suggested by Mumford, we would like to
give reference to earlier work
by Zarkhin. Using techniques similar to ours,
he introduced in \cite{za82} the general notion
of a theta structure compatible with a given endomorphism.
His notion of `compatibility' differs from the one that
we introduce in here. As a consequence, the
resulting formulae (see \cite[Prop.3.2.4]{za82}
and \cite[Th.4.3.1]{za82}) are
significantly different from ours (see Corollary \ref{rmsqrt3eqs}).
An important difference between his and our system
of equations is given by the larger number of coordinates of the projective
space in that the corresponding Hilbert-Blumenthal
variety is embedded. We claim that our formulae are more suitable
for algorithmic applications.
In fact, using our formulas
we are able to compute large explicit examples (see
Section \ref{exx}).
\newline\indent
The problem of computing defining equations for abelian varieties
with real multiplication can also be addressed by methods
using invariants of hyperelliptic genus $2$ curves.
For example, in recent work \cite{grue} Gruenewald
managed to compute the defining equations of the moduli
space of hyperelliptic curves with real multiplication by
$\sqrt{d}$ in terms of Rosenhain invariants for $d \leq 17$.
His computations are based on Runge's method \cite{ru99}.
Gruenewald's approach is as follows.
One considers a suitable system of functions on the
moduli space of hyperelliptic genus-$2$ curves.
Assuming that a power series expansion of these
functions is known, one tries to find linear
relations of given degree
between monomials in the truncated power series
by solving a linear system.
Since the coefficients of the relations grow fast with the discriminant
and his method is exponential in the total degree,
Gruenewald can only compute
real multiplication equations for relatively small discriminant.
Our hope is that a potential generalization of our method for $\sqrt{d}$
multiplication would lead to a more efficient method
which is effective for higher discriminants.
\newline\indent
In the following we explain our strategy
for the computation of equations defining
Hilbert-Blumenthal varieties, i.e. the
moduli spaces of abelian surfaces with specific real multiplication
type. Consider the moduli space $\mathcal{A}_{2,4}^{\Theta}$
of abelian surfaces with level-$4$ theta structure
(compare \cite[$\S$6]{mu66} and \cite[$\S$4]{kempf89}).
A quasi-projective model of $\mathcal{A}_{2,4}^{\Theta}$ can be
given as open immersion into an intersection of quartics. These equations
are given below, see equations (\ref{mum}).
In this article we compute a $\sqrt{3}$-correspondence,
which forms a subspace of the product space $\mathcal{A}_{2,4}^{\Theta}
\times \mathcal{A}_{2,4}^{\Theta}$. The latter
subspace consists of all pairs
of points such that the two points belong to abelian surfaces
with intermediate isogeny of $\sqrt{3}$-type.
Under suitable assumptions regarding the action of the
$\sqrt{3}$-isogeny on line bundles, theta structures and
induced Lagrangian level
structures, one obtains an additional symmetry,
in form of linear relations,
which enable one to turn the correspondence into
a subspace of $\mathcal{A}_{2,4}^{\Theta}$ given
by equations of type (\ref{corr}).
In order to compute the $\sqrt{3}$-correspondence
we apply the multiplication formula for $3$-tuples
of algebraic theta functions (compare \cite[$\S$3.5]{ckl08}).
Such a multiplication formula exists for arbitrary
products of theta functions.
This justifies our expectation that by means of
the general multiplication formula one can compute defining
equations for $\sqrt{d}$-multiplication,
where $d \geq 5$ is an odd square-free integer.
Also, we expect that our method
has some generalization to higher dimension. 
\newline\indent
As we will see in the following, our equations have trivial
coefficients, are of particular simple
shape and bear a striking symmetry.
The latter is induced by elementary combinatorial data.
In the following we give the details of our results.
Let us first fix some notation.
By $A_{2,4}^{\Theta}$ we denote
the moduli space of triples
$(A,\pol,\Theta)$ where
\begin{enumerate}
\item
$A$ is an abelian surface,
\item
$\pol$ is a relatively ample, normalized and
totally symmetric line bundle with $A[\pol]=A[2]$,
\item
$\Theta$ is a symmetric theta structure of type $(\xz / 4 \xz)^2$
for the line bundle $\pol^2$.
\end{enumerate}
There exist canonical coordinates on the
space $A_{2,4}^{\Theta}$ which are induced by
evaluating a canonical basis of the line
bundle $\pol^2$ at the zero section of $A$ for
every triple $(A,\pol,\Theta)$ as above.
In this way, any such triple gives rise to a point $(a_u)_{u \in (\xz / 4 \xz)^2}$
in $\mathbb{P}^{15}$.
Recall that the space $A_{2,4}^{\Theta}$  forms an open
subscheme of the scheme which is defined by the equations
(compare \cite[$\S$6]{mu67a})
\begin{eqnarray}
\label{mum}
(a_{00}^2+a_{02}^2+a_{20}^2+a_{22}^2)(a_{00}a_{02}+a_{20}a_{22}) 
  & = & 2(a_{01}^2+a_{21}^2)^2 \\
\nonumber (a_{00}^2+a_{02}^2+a_{20}^2+a_{22}^2)(a_{00}a_{20}+a_{02}a_{22}) 
  & = & 2(a_{10}^2+a_{12}^2)^2 \\
\nonumber (a_{00}^2+a_{02}^2+a_{20}^2+a_{22}^2)(a_{00}a_{22}+a_{20}a_{02}) 
  & = & 2(a_{11}^2+a_{13}^2)^2 \\
\nonumber (a_{00}a_{20}+a_{02}a_{22})(a_{00}a_{22}+a_{02}a_{20}) 
  & = & 4a_{01}^2 a_{21}^2 \\
\nonumber (a_{00}a_{02}+a_{20}a_{22})(a_{00}a_{22}+a_{02}a_{20}) 
  & = & 4a_{10}^2 a_{12}^2 \\
\nonumber (a_{00}a_{02}+a_{20}a_{22})(a_{00}a_{20}+a_{02}a_{22}) 
  & = & 4a_{11}^2 a_{13}^2 \\
\nonumber (a_{00}^2+a_{02}^2+a_{20}^2+a_{22}^2)a_{13}a_{11}
  & = & (a_{12}^2+a_{10}^2)(a_{01}^2+a_{21}^2) \\
\nonumber (a_{00}^2+a_{02}^2+a_{20}^2+a_{22}^2)a_{01}a_{21} 
  & = & (a_{12}^2+a_{10}^2)(a_{11}^2+a_{13}^2) \\
\nonumber (a_{00}^2+a_{02}^2+a_{20}^2+a_{22}^2)a_{10}a_{12}
  & = & (a_{01}^2+a_{21}^2)(a_{11}^2+a_{13}^2) \\
\nonumber (a_{02}a_{20}+a_{00}a_{22})a_{11}a_{13}
  & = & 2 a_{01}a_{10}a_{21}a_{12} \\
\nonumber (a_{20}a_{00}+a_{22}a_{02})a_{10}a_{12} 
  & = & 2a_{11}a_{13}a_{21}a_{01} \\
\nonumber (a_{00}a_{02}+a_{20}a_{22})a_{21}a_{01} 
  & = & 2a_{11}a_{13}a_{10}a_{12} \\
\nonumber (a_{02}a_{20}+a_{00}a_{22})(a_{01}^2+a_{21}^2)
  & = & 2 a_{10}a_{12}(a_{11}^2+a_{13}^2) \\
\nonumber (a_{00}a_{02}+a_{20}a_{22})(a_{11}^2+a_{13}^2)
  & = & 2a_{10}a_{12}(a_{01}^2+a_{21}^2) \\
\nonumber (a_{02}a_{20}+a_{00}a_{22})(a_{10}^2+a_{12}^2)
  & = & 2 a_{21}a_{01}(a_{11}^2+a_{13}^2) \\
\nonumber (a_{20}a_{00}+a_{22}a_{02})(a_{13}^2+a_{11}^2)
  & = & 2a_{21}a_{01}(a_{10}^2+a_{12}^2) \\
\nonumber (a_{20}a_{00}+a_{22}a_{02})(a_{21}^2+a_{01}^2)
  & = & 2a_{11}a_{13}(a_{10}^2+a_{12}^2) \\
\nonumber (a_{00}a_{02}+a_{20}a_{22})(a_{12}^2+a_{10}^2)
  & = & 2a_{11}a_{13}(a_{01}^2+a_{21}^2) \\
\nonumber a_{01}a_{21}(a_{01}^2+a_{21}^2) & = & a_{10}a_{12}(a_{10}^2+a_{12}^2) \\
\nonumber a_{01}a_{21}(a_{01}^2+a_{21}^2) & = & a_{11}a_{13}(a_{11}^2+a_{13}^2) \\
\nonumber a_{11}=a_{33}, && a_{10}=a_{30} \\
\nonumber a_{01}=a_{03}, && a_{13}=a_{31} \\
\nonumber a_{32}=a_{12}, && a_{21}=a_{23}
\end{eqnarray}
In the following we give the moduli space of abelian surfaces with
real multiplication by $\sqrt{3}$ as a subspace
of the space which is cut out be the relations (\ref{mum}).
We say that a triple $(A,\pol,\Theta)$ is $\sqrt{3}$-admissible
if
\begin{enumerate}
\item
$A$ possesses an endomorphism $\sqrt{3}:A \rightarrow A$ with
$\sqrt{3} \circ \sqrt{3} = [3]$,
\item
$\sqrt{3}^* \pol \cong \pol^3$,
\item
$\sqrt{3}$ acts neutral diagonally on $\Theta$
by the matrix $\squaremat{0}{3}{1}{0}$. 
\end{enumerate} 
We prove in this article that,
if $(A,\pol,\Theta)$ is a $\sqrt{3}$-admissible triple, then the
corresponding theta null point $(a_u)_{u \in (\xz / 4 \xz)^2}$
satisfies the relations
\begin{eqnarray}
\label{corr}
&& 0= a_{13}^2-a_{10}a_{21}-a_{12}a_{01}+a_{11}^2 \\
\nonumber && 0= -a_{00}a_{02}+2a_{10}a_{01}-a_{22}a_{20}-a_{02}a_{22} -a_{20}a_{00}+2a_{12}a_{21}\\
\nonumber && 0= -a_{20}^2-2a_{00}a_{22}-a_{02}^2+4a_{13}a_{11}
\end{eqnarray}
Combining the equations (\ref{corr}) and (\ref{mum}) we obtain
a scheme of relative dimension $2$ over $\xz[1/6]$ whose points correspond
to abelian surfaces with $\sqrt{3}$-multiplication. Our
moduli space forms a finite cover of the above Shimura variety
$\mathrm{Sh}(G,X)_H$ where $G$, $X$ and $H$ are as above.

\subsection{Leitfaden}

This article is structured as follows.
For general type we give the square root three equations
in Theorem \ref{sqrt3mult} of Section \ref{correspond}.
A notion of $\sqrt{d}$-admissibility for theta structures is developed in
the Sections \ref{genrem}, \ref{diagauto} and \ref{admtheta}.
In Section \ref{squares} we describe the locus of split abelian
surfaces.
Finally, in the Sections \ref{goingup} and \ref{thomae} we recall some
useful formulas which give a computable version of
the Torelli morphism. Using these formulas we were able to
compute the examples of curves with $\sqrt{3}$-multiplication,
which are given in Section \ref{exx}.

\subsection{Notation}

A finite faithfully flat morphism of abelian schemes
is called an \emph{isogeny}.
\newline\indent
Let $A$ be an abelian scheme.
For a natural number $m \geq 1$ we denote the multiplication-by-$m$
morphism on $A$
by $[m]_A$, or simply by $[m]$.
Let $\mathcal{L}$ be a line bundle on $A$ and
let $\varphi_{\mathcal{L}}:A \rightarrow \mathrm{Pic}^0_A$
be given by $x \mapsto \langle t_x^* \mathcal{L} \otimes \mathcal{L}^{-1}
\rangle$.
We denote $A[\mathcal{L}]=\mathrm{ker}(\varphi_{\mathcal{L}})$
and $A^{\dagger}=\mathrm{Pic}^0_A$.
An isogeny $\varphi:A \rightarrow A^{\dagger}$ is called
a polarization if locally for the fppf-topology the
morphism $\varphi$ is of the form $\varphi_{\pol}$
for some relatively ample line bundle $\pol$.
A polarization $\varphi$ is called principal if
$\varphi$ gives an isomorphism $A \stackrel{\sim}{\rightarrow} A^{\dagger}$.
For a morphism
$I:A_1 \rightarrow A_2$ we denote the dual morphism $\mathrm{Pic}^0(I):\mathrm{Pic}^0_{A_2} \rightarrow \mathrm{Pic}^0_{A_1}$
by $I^{\dagger}$.
\newline\indent
We set $Z_n=(\xz/n \xz)^2$ for all integers $n \geq 1$.
We define the functor of finite theta functions $V \big( Z_n \big)$ of
type $Z_n$ over a ring $R$ as the sheaf of morphisms
$\underline{\mathrm{Mor}} \big( Z_{n,R}
,\mathcal{O}_R \big)$.
We denote the Cartier dual of the constant group $Z_{n,R}$
by $\hat{Z}_{n,R}$.

\section{Pull back of line bundles under real multiplications}
\label{genrem}

Let $A$
be an abelian surface over a local ring $R$, and let
$\mathcal{L}$ be a normalized and
ample line bundle of degree $1$ on $A$.
Assume that we are given an endomorphism $\sqrt{d}:A \rightarrow A$
such that
\begin{eqnarray}
\label{idem}
\sqrt{d} \circ \sqrt{d} = [d]
\end{eqnarray}
for a given odd square-free integer $d \geq 3$.
In this case we say that $A$ has a $\sqrt{d}$-multiplication.
For simplicity, we suppose that $2d$
is invertible in $R$.
We set $\mathcal{M}=\pol \otimes [-1]^* \pol$.
We remark that the line bundle $\mathcal{M}$ is totally
symmetric. Recall that a line bundle is called
totally symmetric if it is the pull back of some line
bundle on the Kummer surface of $A$. 
\begin{lemma}
\label{pullback}
Possibly after an \'etale extension of the base, one can choose an
ample line bundle
$\mathcal{L}_{\sqrt{d}}$ on $A$ such that
for $\mathcal{M}_{\sqrt{d}}=\mathcal{L}_{\sqrt{d}} \otimes [-1]^* \mathcal{L}_{\sqrt{d}}$
we have
\[
\big( \sqrt{d} \big)^* \mathcal{M}_{\sqrt{d}} \cong \mathcal{M}^{d}
\quad \mbox{and} \quad \big( \sqrt{d} \big)^* \mathcal{M} \cong \mathcal{M}_{\sqrt{d}}^d
\]
\end{lemma}
\begin{proof}
Consider the line bundle $\pol^d$. We have $A[\pol^d]=A[d]$.
Let $K_{\sqrt{d}}$ denote the kernel of the isogeny $\sqrt{d}$. 
In order to descend the line bundle $\pol^d$ along the
isogeny $\sqrt{d}$ one has to construct a lift $\tilde{K}_{\sqrt{d}}$
of $K_{\sqrt{d}}$ to the theta group $G(\pol^d)$
(compare \cite[$\S$1,Prop.1]{mu66}). Since by assumption $2d$ is invertible in $R$,
one knows that such a lift $K_{\sqrt{d}}$ exists over a local \'etale extension
of $R$.
The lift $K_{\sqrt{d}}$ induces a line bundle $\pol_{\sqrt{d}}$
on $A$ such that $\big( \sqrt{d} \big)^* \pol_{\sqrt{d}} \cong \pol^d$.
We note that the possible choices for $\tilde{K}_{\sqrt{d}}$ form a torsor under
the Cartier dual $\hat{K}_{\sqrt{d}}$ of $K_{\sqrt{d}}$.
Hence also the possible line bundles $\pol_{\sqrt{d}}$ form
a $\hat{K}_{\sqrt{d}}$-torsor, so that
there is no canonical choice for the line bundle $\pol_{\sqrt{d}}$.
As an immediate consequence of the above discussion,
we deduce the first part of our claim
\[
\big( \sqrt{d} \big)^* \mathcal{M}_{\sqrt{d}} \cong \mathcal{M}^d.
\]
Furthermore, we claim that
$\big( \sqrt{d} \big)^* \mathcal{M} \cong \mathcal{M}_{\sqrt{d}}^d$.
We know by the symmetry of $\mathcal{M}$ that
\begin{eqnarray}
\label{prop1}
\big( \sqrt{d} \big)^* \left( \big(\sqrt{d} \big)^* \mathcal{M}
\right) \cong [d]^* \mathcal{M} \cong
\mathcal{M}^{d^2} \cong \left( \big( \sqrt{d} \big)^*\mathcal{M}_{\sqrt{d}} \right)^{d}
\cong \big( \sqrt{d} \big)^*\mathcal{M}_{\sqrt{d}}^{d}
\end{eqnarray}
The equality (\ref{prop1}) tells us that the line bundle
$\mathcal{N}= \big( \sqrt{d} \big)^* \mathcal{M}
\otimes \mathcal{M}_{\sqrt{d}}^{-d}$ gives an element
of $\mathrm{ker} \left( \big( \sqrt{d} \big)^{\dagger} \right)
\subseteq A^{\dagger}[d]$.
Also, the symmetry of the line bundles
$\big( \sqrt{d} \big)^* \mathcal{M}_{\sqrt{d}}$ and $\mathcal{M}^{d}$
implies that the class of $\mathcal{N}$ is contained in $A^{\dagger}[2]$.
We conclude from $(d,2)=1$ that
$\mathcal{N}$ induces the trivial class in $\mathrm{Pic}^0_A$.
This proves the lemma.
\end{proof}
\noindent
The following example shows
that in general one has $\mathcal{M} \not\cong \mathcal{M}_{\sqrt{d}}$
for an ample totally symmetric line bundle $\mathcal{M}$.
\begin{example}
Let $E$ be a non-singular
elliptic curve over an algebraically closed field $k$,
where $\mathrm{char}(k) >3$.
We define an isogeny $\sqrt{3}:E^2 \rightarrow E^2$ by
the formula $(x,y) \mapsto (x+2y,x-y)$.
It has the property that $[3]=\sqrt{3} \circ \sqrt{3}$, where
$[3]$ denotes the multiplication-by-$3$ isogeny on $E^2$.
The kernel is given by the map
\[
E[3] \hookrightarrow
E^2, \hspace{0.1cm} x \mapsto (x,x).
\]
Consider the ample symmetric line bundle $\mathcal{L}$ on $E$
which corresponds to the divisor
$2(0_E)$, where $0_E$ denotes the zero section of $E$.
Let $p_i:E^2 \rightarrow E$ $(i=1,2)$ denote the projection
on the $i$-th factor.
We set $\mathcal{M}=p_1^* \mathcal{L} \otimes p_2^* \mathcal{L}$.
The line bundle $\mathcal{M}$ is ample and totally symmetric
with $E^2[\mathcal{M}] = E^2[2]$.
We claim that $ \big( \sqrt{3} \big)^* \mathcal{M}
\not\cong \mathcal{M}^3$. 
Our claim follows from the following computation.
We define a morphism $s:E \rightarrow E^2$ by setting
$x \mapsto (x,0)$.
Then we have
\begin{eqnarray*}
\lefteqn{s^* \mathcal{M}^3 \cong \big( s^* \mathcal{M} \big)^3
\cong \big( s^*(p_1^* \pol \otimes p_2^* \pol)   \big)^3} \\
& & \big( (p_1 \circ s)^* \pol \otimes (p_2 \circ s)^* \pol \big)^3
\cong \pol^3
\end{eqnarray*}
and
\begin{eqnarray*}
\lefteqn{
s^* \left( \big( \sqrt{3} \big)^* \mathcal{M} \right)
\cong s^* \left( \big( \sqrt{3} \big)^* \left( p_1^* \pol  \right)
\right) \otimes
s^* \left( \big( \sqrt{3} \big)^* \left( p_2^* \pol
\right) \right) } \\
& &
\cong (p_1 \circ \sqrt{3} \circ s)^* \pol
\otimes (p_2 \circ \sqrt{3} \circ s)^* \pol
\cong \pol \otimes \pol \cong \pol^2.
\end{eqnarray*}
This proves that
$\mathcal{M} \not\cong \mathcal{M}_{\sqrt{3}}$.
\end{example}
\noindent
Let $R$ and $d$ be as above.
Now suppose that we are given a polarized
abelian scheme $A$ over $R$ with $\sqrt{d}$-multiplication.
Denote the polarization by $\varphi:A \rightarrow A^{\dagger}$.
\begin{definition}
We say that $\varphi$ is $\sqrt{d}$-linear if
$\varphi \circ \sqrt{d} = \big( \sqrt{d} \big)^{\dagger} \circ \varphi$.
\end{definition}
\noindent
Suppose that we are given a
principal polarization $\varphi:A \stackrel{\sim}{\rightarrow} A^{\dagger}$
such that $\varphi=\varphi_{\pol}$ for some normalized ample
line bundle $\pol$ on $A$.
Let $\mathcal{M}$ and $\mathcal{M}_{\sqrt{d}}$ be as in Lemma \ref{pullback}.
\begin{lemma}
The polarization $\varphi_{\pol}$ is $\sqrt{d}$-linear if
and only if $\mathcal{M} \cong \mathcal{M}_{\sqrt{d}}$.
\end{lemma}
\begin{proof}
First, we note that because $[-1]^* \pol$ is algebraically
equivalent to $\pol$ the line bundle $\mathcal{M}$
induces the polarization $\varphi_{\mathcal{M}}=[2] \circ \varphi_{\pol}$.
We conclude that $\varphi_{\pol}$ is $\sqrt{d}$-linear
if and only if $\varphi_{\mathcal{M}}$ is $\sqrt{d}$-linear.
\newline\indent
Suppose that $\varphi_{\mathcal{M}}$ is $\sqrt{d}$-linear. Be definition we have
\begin{eqnarray}
\label{mlin}
\varphi_{\mathcal{M}} \circ \sqrt{d} = \big( \sqrt{d} \big)^{\dagger} \circ \varphi_{\mathcal{M}}
\end{eqnarray}
Let $y \in A$ and set $x=\sqrt{d}(y)$.
We denote the translation-by-$x$ and translation-by-$y$ with the symbols $t_x$ and $t_y$,
respectively.
Then the equation (\ref{mlin}) is equivalent to
the equality
\begin{eqnarray*}
\lefteqn{ \varphi_{\mathcal{M}^d}(y)=\big( \varphi_{\mathcal{M}} \circ [d] \big)(y)=
\big( \varphi_{\mathcal{M}} \circ \sqrt{d} \big) (x)
\stackrel{\rm (\ref{mlin})}{=} \Big(
\big( \sqrt{d} \big)^{\dagger} \circ \varphi_{\mathcal{M}} \Big) (x)} \\
&& = \big( \sqrt{d}\big)^* \big( \langle t_x^* \mathcal{M} \otimes \mathcal{M}^{-1} \rangle \big)
= \langle \big( \sqrt{d} \big)^*t_x^* \mathcal{M} \otimes \big(\sqrt{d}\big)^*\mathcal{M}^{-1} \rangle \\
&&
= \langle (t_x \circ \sqrt{d})^* \mathcal{M} \otimes \big( \sqrt{d} \big)^* \mathcal{M}^{-1} \rangle
= \langle (\sqrt{d} \circ t_y )^* \mathcal{M} \otimes \big( \sqrt{d} \big)^* \mathcal{M}^{-1} \rangle \\
&& = \langle t_y^* \big( \sqrt{d} \big)^* \mathcal{M} \otimes \big(\sqrt{d}\big)^*\mathcal{M}^{-1} \rangle = \varphi_{( \sqrt{d})^* \mathcal{M}}(y).
\end{eqnarray*}
The claim now follows from the fact that two ample totally symmetric
line bundles which give the same polarization are in fact isomorphic. 
\end{proof}
\noindent
As above, let $R$ be a local ring and $d \geq 3$ an odd square-free integer
such that $2d$ is invertible in $R$. 
Assume that we are given a pair $(A,\mathcal{M})$ where
\begin{enumerate}
\item
$A$ is an abelian scheme over the ring $R$, together with
an endomorphism $\sqrt{d}:A \rightarrow A$ such that
$\sqrt{d} \circ \sqrt{d}=[d]$,
\item
$\mathcal{M}$ is a normalized, ample and totally symmetric
line bundle on $A$ such that $A[\mathcal{M}]=A[2]$.
\end{enumerate}
\begin{definition}
The pair $(A,\mathcal{M})$ is called $\sqrt{d}$-admissible
if $\big( \sqrt{d} \big)^* \mathcal{M} \cong \mathcal{M}^d$.
\end{definition}
As we will see in the following sections,
the notion of a $\sqrt{d}$-admissible line bundle is essential for
the construction of a universal abelian surface
with real multiplication by $\sqrt{d}$. 

\section{Diagonalized automorphisms of the theta group}
\label{diagauto}

In the following let $R$ be a local ring such that $2 \in R^*$.
Assume that we are given an abelian scheme $A$ over a local
ring $R$ endowed with a normalized, ample and totally symmetric
line bundle $\mathcal{M}$ such that $A[\mathcal{M}]=A[2]$.
Let $n \geq 1$ be an integer and
let $\Theta_1$ and $\Theta_2$ be theta structures of type $Z_{2n}$
for the line bundle $\mathcal{M}^n$
and let $\tau \in \mathrm{Aut} \big( G(Z_{2n}) \big)$
be a $\xg_{m,R}$-equivariant isomorphism
such that
\begin{eqnarray}
\label{linked}
\Theta_2= \Theta_1 \circ \tau.
\end{eqnarray}
In the following we denote $H(Z_{2n})=Z_{2n} \times \hat{Z}_{2n}$.
The morphism $\tau$ decomposes as a triple
$(\tau_0,\tau_1,\tau_2)$, where
the first component $\tau_0$ satisfies
\begin{eqnarray}
\label{groepswet}
\lefteqn{\tau_0(1,x_1+x_2,l_1l_2) } \\
\nonumber && =e_2 \big( (\tau_1(1,x_1,l_1), l_2),
(x_1,\tau_2(1,x_2,l_2)) \big) \cdot \tau_0(1,x_1,l_1) \cdot \tau_0(1,x_2,l_2)
\end{eqnarray}
for all $x_1,x_2 \in Z_{2n}$ and $l_1,l_2 \in \hat{Z}_{2n}$,
the second component $\tau_1$ is a homomorphism $H(Z_{2n}) \rightarrow Z_{2n}$ and
the third component $\tau_2$ is a homomorphism $H(Z_{2n}) \rightarrow \hat{Z}_{2n}$,
such that $\tau_1 \times \tau_2$ is compatible with
commutator pairings.
\begin{definition}
\label{ndiag}
We call an automorphism $\tau$ as above a diagonalized automorphism if
if there
exists a $\psi \in \mathrm{Aut}(Z_{2n})$ such that
$\tau_1(1,x,l)= \psi(x)$ and $\tau_2(1,x,l)= \hat{\psi}^{-1}(l)$,
where $\hat{\psi}:\hat{Z}_{2n} \rightarrow \hat{Z}_{2n}$
denotes the dual morphism of $\psi$.
In the latter case we denote $\tau=\tau(\psi)$.
Furthermore, if $\tau$ is diagonalized and $\tau_0$ is trivial, then we call
$\tau$ a neutral diagonal automorphism.
\end{definition}
\noindent
We denote the projective theta null points with respect to
the theta structures $\Theta_1$ and $\Theta_2$
by $\big( a_{1,u} \big)_{u \in Z_{2n}}$ 
and $\big( a_{2,u} \big)_{u \in Z_{2n}}$.
\begin{lemma}
\label{dings}
If $\tau \in \mathrm{Aut} \big( G(Z_{2n}) \big)$ is a
diagonal automorphism, then
there exists a $\lambda \in R^*$ such that
\[
a_{2,v}^{2n} = \lambda \cdot a_{1,\psi(v)}^{2n}
\]
for all $v \in Z_{2n}$.
Moreover, if $\tau$ is neutral diagonal, then
there exists a $\lambda \in R^*$ such that
\[
a_{2,v} = \lambda \cdot a_{1,\psi(v)}
\]
for all $v \in Z_{2n}$.
\end{lemma}
\begin{proof}
Assume that we have chosen theta group equivariant isomorphisms
\[
\mu_i:\pi_* \mathcal{M}^n \stackrel{\sim}{\rightarrow} V(Z_{2n}),
\quad i=1,2 
\]
induced by the theta structures $\Theta_1$ and $\Theta_2$,
respectively. Here $\pi$ denotes the structure morphism of $A$.
\newline\indent
The natural action $\star$ of $G \big( \mathcal{M}^n \big)$ on
$\pi_* \mathcal{M}^n$ carries over via $\mu_i$
to the standard action
\[
\big( (\alpha,x,l) \star f \big)(y)
= \alpha \cdot l(y) \cdot f(x+y)
\]
of $G(Z_{2n})$ on $V(Z_{2n})$. 
Let $(\delta_u)_{u \in Z_{2n}}$ denote the Dirac basis
of $V(Z_{2n})$ and $(s^{(i)}_u)_{u \in Z_{2n}}$
the basis of $\pi_* \mathcal{M}^n$ induced by $\mu_i$.
Now let $K=\Theta_2(1,0,\hat{Z}_{2n})$.
We set $t_0 = \sum_{g \in K} g \star s_0^{(1)}$. Then
$t_0$ is invariant under the action of $(1,0,\hat{Z}_{2n})$ via
the theta structure $\Theta_2$.
Let $t_x$ denote the section $\Theta_2(1,-x,1) \star t_0$.
Then $(t_x)_{x \in Z_{2n}}$ forms a basis of $\pi_* \mathcal{M}^n$ 
which induces a theta group equivariant
isomorphism $\mu_2':\pi_* \mathcal{M}^n
\stackrel{\sim}{\rightarrow} V(Z_{2n})$
by mapping $t_x \mapsto \delta_x$.
By uniqueness we conclude that $\mu_2$ and $\mu_2'$
differ by a unit. As a consequence, the
sections $t_0$ and $s_0^{(2)}$ differ by a unit.
In fact, by renormalizing $\mu_2$ one can assume that $t_0=s_0^{(2)}$.
We calculate
\begin{eqnarray*}
\lefteqn{s_x^{(2)}=\Theta_2(1,-x,1) \star t_0} \\
&&  = \Theta_2(1,-x,1) \star \left( \sum_{g \in K} g \star s_0^{(1)} \right)
= \sum_{l \in \hat{Z}_{2n}} \Theta_2(l(-x),-x,l) \star s_0^{(1)} \\
&& \stackrel{{\rm (\ref{linked})}}{=} \sum_{l \in \hat{Z}_{2n}} 
(\Theta_1 \circ \tau) \big( l(-x),-x,l \big) \star s_0^{(1)} \\
&& = \sum_{l \in \hat{Z}_{2n}} 
\Theta_1 \big( l(-x) \cdot \tau_0(1,-x,l), \tau_1(1,-x,l),\tau_2(1,-x,l) \big) \star s_0^{(1)}.
\end{eqnarray*}
This implies
\begin{eqnarray}
\label{basiceq}
\nonumber \lefteqn{ \mu_1 \big( s_x^{(2)} \big) = \mu_1 \left( \sum_{l \in \hat{Z}_{2n}} 
\Theta_1 \big( l(-x) \cdot \tau_0(1,-x,l), \tau_1(1,-x,l),\tau_2(1,-x,l) \big) \star s_0^{(1)} \right) } \\
\nonumber && = \sum_{l \in \hat{Z}_{2n}} 
\big( l(-x) \cdot \tau_0(1,-x,l), \tau_1(1,-x,l),\tau_2(1,-x,l) \big) \star \mu_1 \big( s_0^{(1)} \big) \\
\nonumber && = \sum_{l \in \hat{Z}_{2n}} 
\big( l(-x) \cdot \tau_0(1,-x,l), \tau_1(1,-x,l),\tau_2(1,-x,l) \big) \star \delta_0 \\
\nonumber && = \sum_{l \in \hat{Z}_{2n}} 
l(-x) \cdot \tau_0(1,-x,l) \cdot \tau_2(1,-x,l) \big( \tau_1(1,-x,l) \big)^{-1} \cdot \delta_{-\tau_1(1,-x,l)} \\
\nonumber && = \sum_{l \in \hat{Z}_{2n}} 
\tau_0(1,-x,l) \cdot e_{2n} \big( (-x,\tau_2(1,-x,l)),(\tau_1(1,-x,l),l) \big) \cdot \mu_1 \big( s^{(1)}_{-\tau_1(1,-x,l)} \big) \\
&& \stackrel{\mathrm{(\ref{groepswet})}}{=} \sum_{l \in \hat{Z}_{2n}} 
\frac{\tau_0(1,-x,l)^3}{\tau_0(1,-2x,l^2)} \cdot \mu_1 \big( s^{(1)}_{-\tau_1(1,-x,l)} \big)
\end{eqnarray}
where $e_{2n}$ denotes the natural commutator pairing on $H(Z_{2n})$ which is given by $e_{2n} \big( (x_1,l_1),(x_2,l_2) \big)= \frac{l_2(x_1)}{l_1(x_2)}$ for
$x_1,x_2 \in Z_{2n}$ and $l_1,l_2 \in \hat{Z}_{2n}$.
\newline\indent
By assumption we have
$\tau_1(1,x,l)=\psi(x)$ and $\tau_2(1,x,l)=\hat{\psi}^{-1}(l)$ for all $(x,l) \in H(Z_{2n})$, where $\psi \in \mathrm{Aut}(Z_{2n})$.
One computes
\begin{eqnarray*}
\lefteqn{
e_{2n} \Big( \big( \tau_1(1,-x,1), l \big),
\big( -x, \tau_2(1,0,l) \big) \Big)} \\
&& = e_{2n} \Big( \big( \psi(-x),l \big),
\big( -x, \hat{\psi}(l) \big) \Big) \\
&& = \frac{ \hat{\psi}^{-1}(l) \big( \psi(-x) \big) }{l(-x)}
=\frac{ \big(l \circ \psi^{-1} \big) (\psi(-x)) }{l(-x)}=1
\end{eqnarray*}
Together with (\ref{groepswet}) the latter equality implies that
\[
\tau_0(1,-x,l)=\tau_0(1,-x,1) \cdot \tau_0(1,0,l)
\]
We note that $\tau_0(1,-x,1)$ and $\tau_0(1,0,l)$
are characters on $Z_{2n}$ and $\hat{Z}_{2n}$, respectively.
It now follows from equation (\ref{basiceq}) that
\begin{eqnarray*}
s_x^{(2)} = s_{\psi(x)}^{(1)} \cdot \tau_0(1,-x,1) \cdot \sum_{l \in \hat{Z}_{2n}} \tau_0(1,0,l)
\end{eqnarray*}
Because the right hand side of the above equation is non-trivial,
we must have $\tau_0(1,0,l)=1$ for all $l \in \hat{Z}_{2n}$.
This implies the lemma.
\end{proof}

\section{Definition $\sqrt{d}$-admissible theta structures}
\label{admtheta}

Let $d \geq 3$ be an odd integer. 
Assume now that we are given a $\sqrt{d}$-admissible
pair $\big( A,\mathcal{M} \big)$ in Section \ref{genrem}.
Suppose that a symmetric
theta structure $\Theta_4$ of type $Z_4$
for $\mathcal{M}^2$ is given.
\begin{lemma}
\label{compatibility}
Possibly over a local \'etale extension of the base,
there exist a theta structure $\Theta_{4d}$
of type $Z_{4d}$ for the line bundle $\mathcal{M}^{2d}$
and a symmetric theta structure
$\Theta_4^{\sqrt{d}}$
of type $Z_{4}$ for the line bundle $\mathcal{M}^2$
such that
\begin{enumerate}
\item
$\Theta_{4d}$ is $d$-compatible with the theta structure $\Theta_{4}$,
\item
$\Theta_{4d}$ is $\sqrt{d}$-compatible with the theta structure
$\Theta_{4}^{\sqrt{d}}$.
\end{enumerate}
There exists
a canonical choice for the theta structure
$\Theta_4^{\sqrt{d}}$.
\end{lemma}
\begin{proof}
For the definition of compatibility and our notation we refer to 
\cite[Def.2.3]{ca05b} and \cite[Def.2.5]{ca05b}.
The proof of the above lemma is an immediate consequence
of \cite[Prop.2.2]{ca05b} and the fact that after a suitable
base extension theta structures of separable type always exist.
In the following we will see that
there exists a natural canonical choice for $\Theta_4^{\sqrt{d}}$
\newline\indent
In the following we will see that
there exists a natural canonical choice for $\Theta_4^{\sqrt{d}}$.
Let us briefly recall the construction of the
theta structure $\Theta_4^{\sqrt{3}}$.
As usual, let $\epsilon_d$ and $E_d$ be the natural
embeddings (compare \cite[$\S$2.3]{ca05b}).
We set $K=\mathrm{Ker} \big( \sqrt{d} \big)$. By the
existence of the theta structure $\Theta_{4d}$, one
can lift the group $K$ to a group $\tilde{K}$
in the theta group $G(\mathcal{M}^{2d})$.
Let $G(\mathcal{M}^{2d})^*$ denote the centralizer
of $\tilde{K}$ in $G(\mathcal{M}^{2d})$.
by Lemma \ref{compatibility} we can assume that $K$, considered as
a subgroup of $Z_{4d} \times \hat{Z}_{4d}$ via
the theta structure $\Theta_{4d}$, splits as
$K=K_1 \times K_2$ with $K_1 \subseteq Z_{4d}$
and $K_2 \subseteq \hat{Z}_{4d}$.
We put $Z'_{4d}=Z_4 \oplus K_1$ and
$\hat{Z}'_{4d}=\hat{Z}_4 \odot K_2$, where we
consider $Z_4$ and $\hat{Z}_4$ as embedded in
$Z_{4d}$ and $\hat{Z}_{4d}$, respectively.
The compatibility of $\Theta_4$, $\Theta_{4d}$ and $\Theta_4^{\sqrt{3}}$
translates to the fact that there exists a commutative diagram
\begin{eqnarray*}
\xymatrix{
G(\mathcal{M}^2) \ar@{->}[r]^{\epsilon_d} \ar@{<-}[d]^{\Theta_4} &
G(\mathcal{M}^{2d}) \ar@{<-}[d]^{\Theta_{4d}}  &
G(\mathcal{M}^{2d})^* \ar@{->}[r]^{\mathrm{can}} \ar@{_{(}->}[l] & G(\mathcal{M}^{2d})^*/\tilde{K} \ar@{<-}[d]_{\cong} \ar@{->}[r]^{\mathrm{can}} &
G(\mathcal{M}^2) \ar@{<-}[d]_{\Theta_4^{\sqrt{d}}} \\
G(Z_4) \ar@{->}[r]^{E_d} & G(Z_{4d}) & 
\xg_{m} \times Z'_{4d} \times
\hat{Z}'_{4d} \ar@{->}[r]  \ar@{->}[u]^{\cong} \ar@{_{(}->}[l] & \xg_{m} \times (Z'_{4d}/K_1) \times
(\hat{Z}'_{4d} /K_2) \ar@{->}[r] & G(Z_4)
}
\end{eqnarray*}
Restricting to $G(Z_4)$, we obtain
a commutative square of isomorphisms
\begin{eqnarray*}
\xymatrix{ G(\mathcal{M}^2) \ar@{->}[r]^{\cong}_{\mathrm{can}} & G(\mathcal{M}^2) \\
G(Z_4) \ar@{->}[u]^{\Theta_4} \ar@{->}[r]^{\cong}
& G(Z_4) \ar@{->}[u]_{\Theta_4^{\sqrt{d}}}}
\end{eqnarray*}
If one assumes that the lower horizontal isomorphism
equals the identity, then this assumption uniquely determines the
theta structure $\Theta_4^{\sqrt{d}}$.
\end{proof}
\noindent
In the following we introduce the concept
of $\sqrt{d}$-admissible theta structure.
Consider the number field $K=\xq[\sqrt{d}]$ with
ring of integers $O_K = \xz \oplus \xz \cdot m_d$, where
\begin{eqnarray*}
m_d= \left\{ \begin{array}{l@{, \quad }l}
\sqrt{d} & d \equiv 3 \bmod 4 \\
\frac{1+ \sqrt{d}}{2}
& d \equiv 1 \bmod 4
\end{array}
\right.
\end{eqnarray*}
Let $M_d$ be the Matrix in $\mathrm{Mat}(2, \xz)$
which is defined as
\[
\squaremat{0}{d}{1}{0}
\quad \mbox{if} \quad d \equiv 3 \bmod 4
\]
and
\[
\squaremat{-1}{\frac{(d-1)}{2}}{2}{1} 
\quad \mbox{if} \quad d \equiv 1 \bmod 4.
\]
We note that a matrix representation of the
multiplication-by-$\sqrt{d}$ endomorphism of $O_K$
with respect to the
$\xz$-basis $\{ 1, m_d \}$
is given by $M_d$.
This motivates the following definition.
Setting $\tau=\Theta_4^{-1} \circ \Theta_4^{\sqrt{d}}$
results in a $\xg_m$-equivariant isomorphism $\tau$ of
the standard theta group $G(Z_4)$.
\begin{definition}
\label{sqrt3like}
We call the triple $(A,\mathcal{M},\Theta_4)$
a $\sqrt{d}$-admissible triple, if
\begin{enumerate}
\item
$(A,\mathcal{M})$ is a $\sqrt{d}$-admissible pair,
\item
$\tau$ is neutral diagonal and $\tau=\tau(M_d)$
\end{enumerate}
{\rm (Notation as in Definition \ref{ndiag})}.
\end{definition}
\noindent
We remark that it is always possible to choose parameters
such that $\big( \Theta^{\sqrt{d}}_4 \big)^{\sqrt{d}}=\Theta_4$.

\section{Theta relations induced by isogeny}
\label{thetaisog}

Let $R$ be a local ring with $2 \in R^*$.
Suppose  $\pi_A:A \rightarrow \mathrm{Spec}(R)$ and
$\pi_B:B \rightarrow \mathrm{Spec}(R)$ are abelian
schemes.
Let $\emm$ be an ample
line bundle on $B$.
Suppose that we are given two $2$-compatible theta structures $\Sigma_{m} : G(Z_{m}) \stackrel{\sim}{\rightarrow} G(\emm)$ and $\Sigma_{2m} : G(Z_{2m}) \stackrel{\sim}{\rightarrow} G(\emm^2)$ for some $m \geq 1$.
Let $I:A \rightarrow B$ be an isogeny and $\# \mathrm{Ker}(I)=d$, where
$d \geq 1$ is an integer.
We set $\pol=I^* \emm$.
Assume that there exist $2$-compatible theta structures $\Theta_{md}:
G(Z_{md}) \stackrel{\sim}{\rightarrow} G(\pol)$ and $\Theta_{2md}:
G(Z_{2md}) \stackrel{\sim}{\rightarrow} G(\pol^2)$
such that $\Theta_{md}$ is $I$-compatible with $\Sigma_{m}$ and $\Theta_{2md}$ is $I$-compatible with $\Sigma_{2m}$.
\newline\indent
By the latter assumption the kernel of $I$  
decomposes as $K_1 \times K_2$ where $K_1$ and
and $K_2$ are contained in the image of $Z_{2md}$ and
$\hat{Z}_{2md}$ under the Lagrangian decomposition
induced by $\Theta_{2md}$.
We denote the orthogonal complement of $K_1 \times K_2$ in
$Z_{2md}$ with respect to the standard alternating pairing by $K_1^{\bot}$.
Because of the $I$-compatibility there exists a surjective homomorphism
$\sigma:K_1^{\bot} \rightarrow Z_{2m}$ such that
$\Sigma_{2m}=\Theta_{2md}(\sigma)$
(Notation as in \cite[$\S$5.2]{ca05b}).
Let $\sigma'$ be the restriction of $\sigma$ to
$K_1^{\bot} \cap Z_{md}$, where we consider
$Z_{md}$ as embedded into $Z_{2md}$.
One can assume that $\Sigma_{m}=\Theta_{md}(\sigma')$.
\newline\indent
By general theory there exist theta group equivariant isomorphisms
\[
\mu_j:\pi_{A,*} \pol^j \stackrel{\sim}{\rightarrow} V(Z_{jmd})
\quad \mbox{and} \quad
\gamma_j:\pi_{B,*} \emm^j \stackrel{\sim}{\rightarrow} V(Z_{jm}).
\]
Suppose that we have chosen rigidifications of the line bundles
$\pol$ and $\emm$. This defines,
by means of $\mu_j$ and $\gamma_j$, theta functions $q_{ \emm^j } \in V(Z_{jm})$ and $q_{ \pol^{j} } \in V(Z_{jmd})$ which interpolate the
coordinates of the corresponding theta null points.
\begin{proposition}
\label{descentone}
There exists a $\lambda \in R^*$ such that for all $x \in K_1^{\bot}$ one has
\begin{eqnarray}
\label{keyformula}
q_{\emm^2} \big( \sigma(x) \big)= \lambda \cdot \sum_{w \in \mathrm{ker}(\sigma)} 
q_{\pol^{2}}(x-w)
\end{eqnarray}
\end{proposition}
\begin{proof}
By Mumford's \emph{$2$-Multiplication Formula} \cite[$\S$3]{mu66}
there exists a $\lambda \in
R^*$ such that for all $x \in Z_{2m}$ we have
\[
(\mathbbm{1} \star \delta_0)(x)= \lambda \cdot \sum_{ y \in x + Z_m} \delta_0(x-y)
q_{\emm^2}(y) = \lambda \cdot q_{\emm^2}(x).
\]
Here $\mathbbm{1}$ denotes the $R$-function which takes the value $1$
on all of $Z_{m}$ and $\delta_0$ is the theta function which is defined as
follows
\[
\delta_0(z)= \left\{ \begin{array}{l@{, \quad}l}
1 & z=0 \\
0 & z \not= 0
\end{array} \right.
\]
for $z \in Z_{m}$.
The \emph{Isogeny Theorem} \cite[$\S$1,Th.4]{mu66}
implies that there exists a $\lambda \in R^*$
such that for $x \in Z_{2dm}$ we have
\[
F^*( \mathbbm{1} \star \delta_0)(x)=
\left\{
\begin{array}{c@{ \quad }c}
\lambda \cdot q_{\emm^2} \big( \sigma(x) \big) & \text{ if } x \in K_1^{\bot} \\
0 & \text{otherwise.}
\end{array}
\right.
\]
Furthermore, there exist $\lambda_1,\lambda_2 \in R^*$ such that for $x \in Z_{md}$ we have
\[
F^*( \mathbbm{1} )(x)=
\left\{
\begin{array}{c@{ \quad }c}
\lambda_1 & \text{ if } x \in K_1^{\bot} \cap Z_{md}\\
0 & \text{otherwise,}
\end{array}
\right.
\]
and
\[
F^*( \delta_0 )(x)=
\left\{
\begin{array}{c@{ \quad }c}
\lambda_2 \cdot \delta_0 \big( \sigma(x) \big) & \text{ if } x \in K_1^{\bot} \cap Z_{md}\\
0 & \text{otherwise.}
\end{array}
\right.
\]
Again by Mumford's multiplication formula there exists a $\lambda
\in R^*$ such that for all $x \in K_1^{\bot}$ we have
\begin{eqnarray*}
\lefteqn{ \big( F^*( \mathbbm{1}) \star F^*(\delta_0) \big) (x)} \\
& & = \lambda \cdot \sum_{ y \in x + Z_{md}} F^*( \mathbbm{1} )(x+y) \cdot F^*(\delta_0)(x-y) \cdot q_{\pol^{2}}(y) \\
& & = \lambda \cdot \lambda_1 \cdot \lambda_2 \cdot \sum_{ w \in \mathrm{Ker}(\sigma)}  q_{\pol^{2}}(x-w).
\end{eqnarray*}
The latter equality holds because $F^*(\delta_0)(x-y)=1$ if $x-y \in \mathrm{Ker}(\sigma)$
and $F^*(\delta_0)(x-y)=0$ otherwise.
The theorem now follows from the observation that
$F^*( \mathbbm{1}) \star F^*(\delta_0)$ and $F^*( \mathbbm{1} \star \delta_0)$
differ by a unit.
\end{proof}

\section{A $\sqrt{3}$-correspondence}
\label{correspond}

Let $A$ be an abelian surface over a local ring $R$, and
let $\mathcal{M}$ be a normalized,
ample and totally symmetric line bundle on $A$ such that
$A[\mathcal{M}]=A[2]$.
Let $\sqrt{3}:A \rightarrow A$ be an isogeny
such that $\sqrt{3} \circ \sqrt{3} = [3]$.
Furthermore, we suppose that $6 \in R^*$.
In the following we assume that $(A,\mathcal{M})$
gives a $\sqrt{3}$-admissible pair.
Assume that we are given a symmetric theta structure $\Theta_4$,
and let $\Theta_4^{\sqrt{3}}$ be the canonical theta structure
of Lemma \ref{compatibility}.
We define
\begin{eqnarray*}
S= \{ (x,y,z) \in Z_4^3 \mid (x-2y,x+y-z,x+y+z) \in Z_2^3 \}.
\end{eqnarray*}
For $(x_1,y_1, z_1),(x_2,y_2,z_2) \in S$ we denote $(x_1,y_1,z_1) \sim (x_2,y_2,z_2)$
if there exists a permutation matrix $P \in \mathrm{Mat}_{3}({\xz})$ such that
\[
(x_1-2y_1,x_1+y_1-z_1,x_1+y_1+z_1) = (x_2-2y_2,x_2+y_2-z_2,x_2+y_2+z_2) P.
\]
Let $(a_u)_{u \in Z_4}$ and $(a_u^{\sqrt{3}})_{u \in Z_4}$ denote the theta null points associated to the theta structures $\Theta_4$ and $\Theta_4^{\sqrt{3}}$,
respectively.
\begin{theorem}
\label{sqrt3mult}
For all pairs of triples
$(x,y_1,z_1),(x,y_2,z_2) \in S$ which satisfy 
$(x,y_1,z_1) \sim (x,y_2,z_2)$ one has
\begin{eqnarray}
\label{type1}
\sum_{ u \in Z_2} a^{\sqrt{3}}_{y_1+u} a_{z_1+u}
=\sum_{ v \in Z_2} a^{\sqrt{3}}_{y_2+v} a_{z_2+v}.
\end{eqnarray}
\end{theorem}
\begin{proof}
There exists a unique theta structure $\Theta_2$
of type $Z_2$ for $(A, \mathcal{M})$ which is $2$-compatible
with the given theta structure $\Theta_4$ (compare \cite[$\S$2,Rem.1]{mu66}).
One can show that there exist theta structures
$\Theta_{6j}$ ($j=1,2$)
of type $Z_{6j}$ for $\mathcal{M}^{3j}$
which are compatible with $\Theta_{2j}$.
In the following we denote $I= \{ 2,3,6 \}$.
For the following we assume that we have chosen rigidifications
for the line bundles $\mathcal{M}^i$
and theta group invariant isomorphisms
\[
\mu_i: \pi_* \mathcal{M}^i \stackrel{\sim}{\rightarrow} V(Z_i), \quad i \in I,
\]
where $\pi:A \rightarrow \mathrm{Spec}(R)$ denotes the
structure morphism.
Our choice determines theta functions $q_{\mathcal{M}^i} \in V(Z_i)$ which
interpolate the coordinates of the theta null point with respect to
$\Theta_i$.
Let $\{ \delta_w \}_{w \in Z_2}$ denote the Dirac basis 
of the module of finite
theta functions $V(Z_2)$.
Let now $(x_0,y_i,z_i) \in S$ where $i=1,2$ and set
\[
(a_i,b_i,c_i)=(x_0-2y_i,x_0+y_i-z_i,x_0+y_i+z_i).
\]
Suppose that $(x_0,y_1,z_1) \sim (x_0,y_2,z_2)$,
i.e. there exists
a permutation matrix $P \in \mathrm{Mat}_3(\xz)$ such that
\begin{eqnarray}
\label{permut}
(a_1,b_1,c_1) = (a_2,b_2,c_2) P.
\end{eqnarray}
For $j=1,2$ we set
\[
S^{(j)}_{x_0}=\{ (x,y,z) \in S \hspace{0.1cm} |  \hspace{0.1cm}
(x=x_0)
\wedge (x-2y,x+y-z,x+y+z)=(a_j,b_j,c_j) \}.
\]
By the $3$-multiplication formula \cite[Th.3.11]{ckl08}
there exists a $\lambda \in R^*$ such that
\begin{eqnarray}
\label{prodbasis}
\lefteqn{ \big( \delta_{a_j} \star \delta_{b_j} \star \delta_{c_j}
  \big) (x_0)} \\
\nonumber & & = \lambda \sum_{ (x,y,z) \in S^{(j)}_{x_0} } \delta_{a_j}(x-2y) \delta_{b_j}(x+y-z)
\delta_{c_j}(x+y+z) q_{\mathcal{M}^{6}}(y) q_{\mathcal{M}^2}(z) \\
\nonumber & & = \lambda \sum_{t \in Z_2}
q_{\mathcal{M}^{6}}(y_j+t) q_{\mathcal{M}^2}(z_j+t).
\end{eqnarray}
The theta structure $\Theta_4^{\sqrt{3}}$
induces a theta group invariant isomorphism
\[
\mu_4^{\sqrt{3}}: \pi_* \mathcal{M}^2 \stackrel{\sim}{\rightarrow} V(Z_4).
\]
Let $q^{\sqrt{3}}_{\mathcal{M}^2} \in V(Z_4)$
be the finite theta function that interpolates the
theta constants $\big( a_u^{\sqrt{3}} \big)_{u \in Z_4}$.
It follows by Proposition \ref{descentone}
that there exists an $\alpha \in R^*$ such that
\begin{eqnarray}
\label{desc}
\sum_{s \in \mathrm{Ker}(\sqrt{d})} q_{\mathcal{M}^{6}}(z+s)= \alpha q^{\sqrt{3}}_{\mathcal{M}^2}(z)
\end{eqnarray}
for all $z \in Z_4$.
Combining the equations (\ref{prodbasis}) and (\ref{desc})
we conclude
that there exists $\lambda \in R^*$ such that
\begin{eqnarray}
\label{form}
\lefteqn{ \sum_{s \in \mathrm{Ker}(\sqrt{d})} \big( \delta_{a_j} \star \delta_{b_j} \star \delta_{c_j} \big) (x_0-s) } \\
\nonumber & & = \lambda \sum_{t \in Z_2} q_{\mathcal{M}^2}(z_j+t)
\sum_{s \in \mathrm{Ker}(\sqrt{d})} q^{\sqrt{3}}_{\mathcal{M}^2}(y_j+t+s) \\
\nonumber & & = \lambda \alpha \sum_{t \in Z_2} q^{\sqrt{3}}_{\mathcal{M}^2} \big( y_j+t \big) q_{\mathcal{M}^2}(z_j+t).
\end{eqnarray}
The commutativity of the $\star$-product and equality (\ref{permut})
imply that
\begin{eqnarray}
\label{comm}
\big( \delta_{a_1} \star \delta_{b_1} \star \delta_{c_1} \big)(x_0-s) =
\big( \delta_{a_2} \star \delta_{b_2} \star \delta_{c_2} \big)(x_0-s).
\end{eqnarray}
As a consequence of the equalities (\ref{form}) and (\ref{comm}) we have
\begin{eqnarray*}
\sum_{u \in Z_2} q^{\sqrt{3}}_{\mathcal{M}^2} \big( y_1+t \big) q_{\mathcal{M}^2}(z_1+u)
= \sum_{v \in Z_2} q^{\sqrt{3}}_{\mathcal{M}^2} \big( y_2+t \big) q_{\mathcal{M}^2}(z_2+v).
\end{eqnarray*}
This implies the equation (\ref{type1}) and
thus completes the proof of the theorem.
\end{proof}
\noindent
We set $M=\squaremat{0}{3}{1}{0}$.
\begin{corollary}
\label{rmsqrt3eqs}
Assume that $\big( A,\mathcal{M},\Theta \big)$
is a $\sqrt{3}$-admissible triple.
Then for all pairs of triples
$(x,y_1,z_1),(x,y_2,z_2) \in S$ with 
$(x,y_1,z_1) \sim (x,y_2,z_2)$ one has
\begin{eqnarray}
\sum_{ u \in Z_2} a_{M(y_1+u)} a_{z_1+u}
=\sum_{ v \in Z_2} a_{M(y_2+v)} a_{z_2+v}.
\end{eqnarray}
\end{corollary}
\noindent
The above corollary is an immediate consequence
of Lemma \ref{dings} and Theorem \ref{sqrt3mult}.

\section{Product abelian surfaces}
\label{squares}

In this section we give equations for the moduli space
of abelian surfaces with theta structure which
decompose as a product of elliptic curves
with product polarization.
\newline\indent
Let $R$ denote a local ring with $2 \in R^*$.
Assume now that we are given two elliptic curves
$E_1$ and $E_2$ over $R$,
which are endowed with normalized, ample and totally symmetric
line bundles $\pol_1$ and $\pol_2$, respectively.
Let $\Theta_i$ be a symmetric theta structure
of type $\xz / 4 \xz$ for the line bundle $\pol_i^2$ ($i=1,2$).
We denote the theta null points induced by
$\Theta_1$ and $\Theta_2$ by $(x_i)_{i \in \xz / 4 \xz}$
and $(y_j)_{j \in \xz / 4 \xz}$, respectively.
We set $A=E_1 \times E_2$ and denote the projection
$A \rightarrow E_i$ by $p_i$. Let $\Theta=\Theta_1 \times
\Theta_2$ be the product theta structure
for the line bundle $\mathcal{M}^2$,
where $\mathcal{M}=p_1^* \pol_1 \otimes p_2^* \pol_2$.
A proof of the existence of the product theta
structure is given in \cite[$\S$3.4]{ckl08}.
Let $(a_{ij})_{i,j \in \xz / 4 \xz}$ denote the
theta null point corresponding to the triple $(A,\Theta,\mathcal{M})$.
It is well-known that one has
\begin{eqnarray}
\label{building}
a_{ij}=x_i \cdot y_j, \quad i,j \in \xz / 4 \xz.
\end{eqnarray}
The theta null points $(x_i)$ and $(y_j)$ satisfy
the relations
\begin{eqnarray*}
(x_0^2+x_2^2)x_0x_2 = 2x_1^4 \quad \mbox{and} \quad
(y_0^2+y_2^2)y_0y_2 = 2y_1^4
\end{eqnarray*}
Multiplying the above two equations with each other we
obtain by substituting relation (\ref{building}) the equation
\begin{eqnarray*}
4a_{11}^4 & = & (a_{00}^2+a_{02}^2+a_{20}^2+a_{22}^2)a_{00}a_{22}
\end{eqnarray*}
Mumford's set of theta relations (\ref{mum}) plus
the latter relation define a
two dimensional fiberwise irreducible space that
equals the moduli space of products of
elliptic curves with product line bundle
and product theta structure.
\newline\indent
Similarly, one obtains the following equations for the
curve of square elliptic curves in $\xp^{15}$ as
subscheme of $A_{2,4}^{\Theta}$ (compare Section \ref{intro})
\begin{eqnarray*}
& a_{11}a_{22}  =  a_{21}^2 & \\
& a_{10}a_{21}  =  a_{11}a_{20} & \\
& a_{10}a_{22}  =  a_{20}a_{21} & \\
& 2a_{11}^2  =  a_{00}a_{20} + a_{20}a_{22} &\\
& a_{00}a_{21}  =  a_{10}a_{20} &\\
& a_{00}a_{11}  =  a_{10}^2 &\\
& a_{00}a_{22}  =  a_{20}^2 &\\
& a_{13} = a_{31}, \; a_{03} = a_{30} &\\
& a_{23} = a_{32}, \; a_{01} = a_{10} &\\
& a_{11} = a_{33}, \; a_{13} = a_{31} &\\
& a_{01} = a_{03}, \; a_{12} = a_{32} &\\
& a_{11} = a_{13}, \; a_{02} = a_{20} &\\
& a_{12} = a_{21}, \; a_{21} = a_{23} &\\
& a_{31} = a_{33}, \; a_{21} = a_{23} &\\
& a_{01} = a_{03}, \; a_{10} = a_{30} &
\end{eqnarray*}

\section{Theta constants of genus-$2$ curves}
\label{thomae}

In this section we recall some classical formulas due to Thomae
which relate the coefficients
of a genus-$2$ curve with its theta constants.
In the following let $R$ be a local ring with maximal
ideal $\mathfrak{m}$ such that $R/ \mathfrak{m}$ is perfect and $2 \in R^*$.
Let $C$ be a smooth hyperelliptic genus-$2$ curve over $R$, i.e. a proper and smooth relative curve such that its geometric fibers have genus $2$.
By assumption, over a local \'etale extension
one can give a plane $R$-model of $C$ which is of the form
\[
y^2z^4=\prod_{i=1}^6 (x-e_iz)
\]
such that $e_i \in R$ and $e_i \not\equiv e_j \bmod \mathfrak{m}$
for $i \not= j$.
We set $J=J(C)=\mathrm{Pic}^0_C$. The abelian scheme $J$
is called the relative Jacobian of $C$.
It is endowed with a
canonical principal polarization $\varphi:J \stackrel{\sim}{\rightarrow}
\mathrm{Pic}^0_J$. Over some local \'etale extension
one can choose a
normalized, ample and totally symmetric line bundle
$\mathcal{L}$ on $J$ such that $\mathcal{L}$
induces $2 \varphi$.
Extending $R$ further we can choose a theta
structure $\Theta_2$ of type
$(\xz / 2 \xz)^2$ for $\mathcal{L}$.
Let $(b_{00},b_{01},b_{10},b_{11})$ be the associated $2$-theta null point.
Then one has
\begin{eqnarray*}
 \big(b_{00}^2+b_{01}^2+b_{10}^2+b_{11}^2 \big)^2 =(e_1-e_3)\cdot(e_1-e_5)\cdot(e_2-e_4)\cdot(e_3-e_5) \\
\big( b_{00}^2-b_{01}^2+b_{10}^2-b_{11}^2 \big)^2 =(e_1-e_3)\cdot(e_1-e_4)\cdot(e_2-e_5)\cdot(e_3-e_4) \\
\big( b_{00}^2+b_{01}^2-b_{10}^2-b_{11}^2 \big)^2 =(e_1-e_2)\cdot(e_1-e_4)\cdot(e_2-e_4)\cdot(e_3-e_5) \\
\big( b_{00}^2-b_{01}^2-b_{10}^2+b_{11}^2 \big)^2 =(e_1-e_2)\cdot(e_1-e_5)\cdot(e_2-e_5)\cdot(e_3-e_4)
\end{eqnarray*}
In \cite[Ch.IIIa,$\S$8,Th.8.1]{tat2} one finds a complex analytic proof of the Thomae formulae.
We note that our formulas differ from the classical Thomae formulae
by an invertible linear transformation. This change of coordinates relates
classical Jacobi- with Mumford-coordinates (compare \cite[p.353]{mu66}).
\newline\indent
The above formulae enable one to compute the $2$-theta null point
of a curve by taking square roots. We note that not all choices
for the square roots lead to a valid $2$-theta null point. 

\section{Going up from level-$2$ to level-$4$}
\label{goingup}

In this section we explain how to compute the theta null point with
respect to the $4$-th power of a line bundle taking square roots of the
coordinates of the theta null point with respect to the $2$nd power.
\newline\indent
Assume that $A$ is an abelian surface
over a local ring $R$ with $2 \in R^*$. Let $\pol$ be a normalized,
ample and totally symmetric line bundle on $A$ such that $A[\pol]=A[2]$.
Suppose that we are given a theta structures $\Theta_2$ of type
$Z_2$ for $\pol$ and a symmetric theta structure
$\Theta_4$ of type $Z_4$ for $\pol^2$.
Let $(b_u)_{u \in Z_2}$ and $(a_v)_{v \in Z_4}$ denote the theta
null points attached to the theta structures $\Theta_2$ and $\Theta_4$,
respectively.
\newline\indent
It follows by Mumford's $2$-multiplication formula that
assuming a suitable normalization one gets
\begin{eqnarray*}
b_{00}^2 & = & (a_{00}^2 + a_{02}^2 + a_{20}^2 + a_{22}^2 ) \\
b_{01}^2 & = & 2 (a_{00}a_{02}+a_{20}a_{22}) \\
b_{10}^2 & = & 2  (a_{00} a_{20} + a_{02} a_{22}) \\
b_{11}^2 & = & 2  (a_{00} a_{22} + a_{02} a_{20}) \\
b_{01} b_{00} & = & 2  (a_{01}^2 + a_{21}^2 ) \\
b_{11} b_{00} & = & 2  (a_{11}^2 + a_{13}^2) \\
b_{10} b_{00} & = & 2  (a_{12}^2 + a_{10}^2) \\
b_{11} b_{01} & = & 4  a_{10} a_{12} \\
b_{11} b_{10} & = & 4  a_{01} a_{21} \\
b_{01} b_{10} & = & 4  a_{11} a_{13}.
\end{eqnarray*}
By the above formulas one has
\begin{eqnarray*}
b_{00}^2+b_{01}^2+b_{10}^2+b_{11}^2 & = & 
(a_{00}+a_{02}+a_{20}+a_{22})^2 \\
b_{00}^2-b_{01}^2+b_{10}^2-b_{11}^2 & = & 
(a_{00}-a_{02}+a_{20}-a_{22})^2 \\
b_{00}^2+b_{01}^2-b_{10}^2-b_{11}^2 & = & 
(a_{00}+a_{02}-a_{20}-a_{22})^2 \\
b_{00}^2-b_{01}^2-b_{10}^2+b_{11}^2 & = & 
(a_{00}-a_{02}-a_{20}+a_{22})^2.
\end{eqnarray*}
Hence we can compute $a_{00},a_{02},a_{20},a_{22}$ from $b_{00},b_{01},b_{10},b_{11}$ by taking square roots over a local \'etale extension of $R$ and solving a
linear system.
\newline\indent
By the above we also have
\begin{eqnarray*}
b_{01} b_{00} + b_{11} b_{10} & = & 2 (a_{01}+a_{21})^2 \\
b_{11} b_{00} + b_{01} b_{10} & = & 2 (a_{11}+a_{13})^2 \\
b_{10} b_{00} + b_{11} b_{01} & = & 2 (a_{10}+a_{12})^2 \\
b_{01} b_{00} - b_{11} b_{10} & = & 2 (a_{01}-a_{21})^2 \\
b_{11} b_{00} - b_{01} b_{10} & = & 2 (a_{11}-a_{13})^2 \\
b_{10} b_{00} - b_{11} b_{01} & = & 2 (a_{10}-a_{12})^2.
\end{eqnarray*}
This allows us to compute theta null values
$a_{01}, a_{11}, a_{10}, a_{12}, a_{13}, a_{21}$
from $b_{00},b_{01},b_{10},b_{11}$ by taking square roots
over a local \'etale extension of $R$ 
and solving a linear system.

\section{Example}
\label{exx}

The real multiplication equations of this article and
the Magma code \cite{magma} that we used to produce the
example are available from the authors website
\begin{verbatim}
http://www.uni-ulm.de/fileadmin/website_uni_ulm/mawi.inst.100/mitarbeiter/
       carls/code/rm_sqrt3.tar.gz
\end{verbatim}
The computation of the following example took
$217$ seconds CPU time on a
Intel Xeon Quad Core system at 3GHz.
This indicates that the generation of abelian surfaces
with RM using Mumford's theta coordinates
is effective for cryptographic size and far beyond.
\newline\indent
Now let $\xf_p$ be a finite field of
characteristic
\begin{eqnarray*}
\lefteqn{p=2800333449468855427078318918911064957466346050652651706861651359541340210937518} \\
&& 9431935592353737690173937903928210822864900855816077651451334148846540061741866 \\
&& 8815872898112505450998897734215679529296359811615259890665122008091641369773466 \\
&& 0509028721903773570791357488028843181423689774686954830610975001
\end{eqnarray*}
Let $\xf_{p^8}=\xf_p[x]/(f)$ where $f(x)=x^8+3$.
In the following we denote the congruence class of $x$ in $\xf_p[x]/(f)$
by $\bar{x}$.
An instance of a
$4$-theta null point of an abelian surface with real multiplication
by $\sqrt{3}$, which satisfies the system of equations
(\ref{mum}) and (\ref{corr}) of Section \ref{intro},
is given in short form by the coordinate tuple
\[
(a_{00},a_{01},a_{02},a_{10},a_{11},a_{12},a_{13},a_{20},a_{21},a_{22})
\]
over $\xf_{p^8}$, where the coefficients are as follows
{\scriptsize
\begin{eqnarray*}
&& \big( 1, 230948181937011261114627330367112087301879751114820109245770593437711125436930 \\
&& 24333055801661362498793076183085439222824115011892917730186589905795938205273950268 \\
&& 79963131158433136866852406009943021958917819236895830196533453409030757762216691186 \\
&& 819798385050850086312507828078505110095419945195774893862*\bar{x}^7 + 2745534299657 \\
&& 46175020413075825439488675085688030727643698858692271770744217837219135744586959810 \\
&& 00470142871456750521658512921119367899354334714961247802906781748810229815778957343 \\
&& 19604975982843854356675373374698422858659311526597199174146577373222835547705601089 \\
&& 171539924525014040772161382829588918662*\bar{x}^3, \\
&& 22730095486682121352411354483037697388765865045297885060737060997306010682832561352 \\
&& 96177013183944678263718967913756788766960014505736233893318803328425461162570939405 \\
&& 39320137749528223807303466589456391156107719286435225989528098408944470455540420527 \\
&& 2675449894780787034386363145157494436293732835314543, \\
&& 12943582435177557375069620841819500547632515379173641673612362030787488252494743861 \\
&& 76329927805790287594896624898126125957553602072369386300031015607788039975910689497 \\
&& 14656366150036322881777788313681383557684571605773053971309320530233521893754611878 \\
&& 4725270567812271550405989627299891975994000976760491*\bar{x}^5 + 214313827016528773 \\
&& 45876302895537694041525165352157389894230439977566645945148264491377816468035422932 \\
&& 95333303123646243904926283910020005995259598022313611608895250530947399818194934731 \\
&& 81211530192709781728847781133715582533426428824772014214107262091321296202226741924 \\
&& 9686378888345142380465002773386405*\bar{x},
\end{eqnarray*}
\begin{eqnarray*}
&& 58416679319653988064159139182552267932412687847858408568636497196105649630443447074 \\
&& 45435552739432128677337851233904011292048284789158369897670142209178581264951771425 \\
&& 02128412125285686036647834690411570248285781378483257085961546769493270562455507506 \\
&& 184862686955304340933187825192517023624642208927715*\bar{x}^4 + 2074036108879803997 \\
&& 18607531300399924133246498405011152989112587276703131289270309311517901140469376982 \\
&& 67485623449059101570879258758180673195264552671008329479344694377038356403104698794 \\
&& 49441909019012688690272139400861677639536508415573518853955253126554580959524702920 \\
&& 482487769984821819753916065326613, \\
&& 10456441303738139864159331416882053654893857621333829033710754076380241336490063119 \\
&& 71774306586121712615882582724769003158024352443766913344735623281817936920055917352 \\
&& 16814357289435015495425786969619417519089003595877462120399744414606626886274135094 \\
&& 1172271067232157627693477430626707345544217495274416*\bar{x}^5 + 395267184170434499 \\
&& 65359971908410560480694102503442726733415492917218081581180921655960884079860607120 \\
&& 36005950876579295214117965686025703443765834840255483846159514031095945029885364066 \\
&& 38309834311653948996264885723882505781204868486376954091553493762521760869863504514 \\
&& 705601391536572919256101989634749*\bar{x}, \\
&& 16611476180322968071442581764726144123091097176406558201264246122910037973247459782 \\
&& 36573816029557453299396167515135304313077697360226592227392708552384250474343223140 \\
&& 49623271369170834148762887996205204665535084018399899386110354886543580271421020509 \\
&& 4140696591450911430362753717918626164403385400697229*\bar{x}^4 + 627454797736226824 \\
&& 89793102007261218690540905944101312928092344557418903122889454521533255837125370767 \\
&& 00693461671371836690331775384453806396574580910070228875267452089512487236091349858 \\
&& 53187554765650944149775483136156983394360078348730054732734904638868724545468587437 \\
&& 782922961345248502054530133246637, \\
&& 17010938916701657456628379671237635855989318577476338484087237808179944112442185980 \\
&& 26365593386886821396855162025912936451213631430541630977303264448572847351583833461 \\
&& 86735463941779495161977870161007659321316835253424312533873955064855470495789895261 \\
&& 9920062579122661783528072962922424317644309768602836, \\
&& 27233231099944079353158981212108525074151340994462561487756803116640731054752310510 \\
&& 38619554472927702027284029068999546114508756704669140633848492785343928801103191428 \\
&& 37150549526162438935884886890745960532415144824599493330363042873481296716808672609 \\
&& 1886496587650573172920713816820914359766929018772135*\bar{x}^7 + 116063463087096368 \\
&& 72381969282865276147339345569360253231261952842158795380628681580765873073520800122 \\
&& 89321744464862438267412395404445533079773488432271940263977653266642554520869281965 \\
&& 44364243784577666416509483148072194218908504283723847336586205849204344151720045439 \\
&& 5640327592754995759419314136530233*\bar{x}^3,\\
&& 179067787751124024504810272450957725606300219066942878644908399290926820281 \\
&& 65226431988350475772431469552080566666743115215560136685115676987590443 \\
&& 06596915983036081508584137568635793550334449713236483566740653822344677 \\
&& 14133320892437663972662096712011419074047709227028778866571162057685105 \\
&& 9606381332281*\bar{x}^4 \big)
\end{eqnarray*}}
The above $4$-theta null point belongs to the Jacobian of a
hyperelliptic genus two curve. 
By the formulas of the Sections \ref{thomae} and \ref{goingup} we can
compute a Rosenhain model
\[
y^2=x(x-1)(x-\lambda_1)(x-\lambda_2)(x-\lambda_3)
\]
of the corresponding genus two curve.
The Rosenhain invariants $\lambda_1,\lambda_2,\lambda_3$
are defined over $\xf_{p^2}=\xf_p[y]/(y^2+3)$.
We denote by $\bar{y}$ the congruence class of $y$ in $\xf_p[y]/(y^2+3)$.
One has
{\scriptsize
\begin{eqnarray*}
\lefteqn{\lambda_1=2039487070036611647680452203516624853229218571374855718055481034790975991504624} \\
&& 534104207145338037519615824772029873706848572462516119500914588865262777614\\
&& 833026966554825050880289809603320699642554137861590919614968437242591975885\\
&& 051898308588368979209664335473578173441866890653892491704979050695061218*\bar{y}\\
&& + 13277853173719792897335784798454693443469601066388025477751482354877482\\
&& 139738408187446973073737501447673971446185559265083096855257637600523086265\\
&& 531442493139361219915488953651292219487491666872545082814224959766539213151\\
&& 510107822625109301785604946714419704613270162600815701451274338842827391960\\
&& 16481 \\
\lefteqn{\lambda_2=1217655172881207009593495496540383150004377001379559710652840633358227029548031} \\
&& 981781518883378396400179705031742337854470869846048516238775213643500600990\\
&& 788578112514082077197258725505108258298123242432208543852008288027377488795\\
&& 13253141074887002117583313260723860278753693081586181113892813240661228*\bar{y} \\
&& + 5115189664835245734672484575020321035910154551257829940510331083207423918\\
&& 297577670444329243496394240297767071320080557495498448221100957232253205725\\
&& 262406520895401644913738885257859665714973221450591234092142012051463407083\\
&& 5979880545838992856365037944123819772191020095694076685123224515478309500837 \\
\lefteqn{\lambda_3=2080768235846809814580109380593740744600107616660538822331578060941667919162528} \\
&& 742182001188848587142284587193632893489607996037046814495252195949978022245\\
&& 354949090312376942167935480700133063373579878022673243157346264565993054103\\
&& 635774984528075902192111080696936002859717856456821884396795922991005097*\bar{y}\\
&& + 22987669990796198593695173822539528890584973931713921607851615480028777\\
&& 883694099560354746831053313724839080186645310895164214547113748420650781676\\
&& 337732615128374738418336460697156731675791188331591100695353029096687691730\\
&& 735222865449183174382616219294841699858227474899144762005138124478726173052\\
&& 04317
\end{eqnarray*}}

\section{Acknowledgments}

We would like to thank David Gruenewald for useful discussions and
for providing real multiplication equations for discriminant $12$
in form of Magma code.
It enabled us to verify the examples of Section \ref{exx}.
\newline\indent
We thank the Communication and Information Center (KIZ, Ulm) staff members,
namely Harald Däubler and Christian Mosch, for providing
support regarding computations on the local compute cluster.  
\newline\indent
I thank Eyal Goren for pointing out the reference to Zarkhin's work.

\bibliographystyle{plain}

\end{document}